\documentclass[a4paper,11pt]{article}

\usepackage{subfig}
\usepackage{pgfplots}
\usepackage{pgfplotstable}
\usepackage{mathtools}
\usepackage{multicol}
\usepackage{comment}
\usepackage{booktabs}
\pgfplotsset{compat=1.5}
\usepackage{amssymb}
\usepackage{url}
\usepackage{bm}

\usepackage{pdfsync}
\usepackage{float}
\usepackage{tabularx}
\usepackage{enumerate}
\usepackage{array}
\usepackage{xspace}

\usepackage{authblk}

\usepackage{amsmath}
\usepackage{tabularx}
\usepackage{booktabs}
\usepackage{url}
\usepackage{tikz}
\usepackage{pgfplots}
\usetikzlibrary{spy}
\usepackage{color}
\newcommand{\TODO}[1]{{\color{black}#1}}

\newcommand{\x}{\ensuremath{\boldsymbol{x}}}
\newcommand{\mupar}{\ensuremath{\boldsymbol{\mu}}}
\newcommand{\etapar}{\ensuremath{\boldsymbol{\eta}}}

\begin{document}

\title{A non-intrusive approach for the reconstruction of POD modal coefficients through active subspaces}

\author[]{Nicola~Demo\footnote{nicola.demo@sissa.it}}
\author[]{Marco~Tezzele\footnote{marco.tezzele@sissa.it}}
\author[]{Gianluigi~Rozza\footnote{gianluigi.rozza@sissa.it}}

\affil[]{Mathematics Area, mathLab, SISSA, International School of Advanced Studies, via Bonomea 265, I-34136 Trieste, Italy}

\maketitle

\begin{abstract}
Reduced order modeling (ROM) provides an efficient framework to compute
solutions of parametric problems. Basically, it exploits a set of precomputed
high-fidelity solutions --- computed for properly chosen parameters, using a
full-order model --- in order to find the low dimensional space that contains
the solution manifold. Using this space, an approximation of the numerical
solution for new parameters can be computed in real-time response scenario,
thanks to the reduced dimensionality of the problem. In a ROM framework, the
most expensive part from the computational viewpoint is the calculation of the
numerical solutions using the full-order model. Of course, the number of
collected solutions is strictly related to the accuracy of the reduced order
model. In this work, we aim at increasing the precision of the model also for few
input solutions by coupling the proper orthogonal decomposition with
interpolation (PODI) --- a data-driven reduced order method --- with the active
subspace (AS) property, an emerging tool for reduction in parameter space. The
enhanced ROM results in a reduced number of input solutions to reach the
desired accuracy. In this contribution, we present the numerical results
obtained by applying this method to a structural problem and in a fluid dynamics
one.
\end{abstract}

% \keywords{Nonintrusive Model Order Reduction, Active Subspaces, Free Form Deformation, POD with
%  Interpolation, Dynamic Mode Decomposition, Parameter Space Reduction}

%\clearpage

\section{Introduction}
\label{sec:intro}

In a large variety of engineering and computational science fields, reduced
order modeling (ROM) has gained more and more popularity to treat parametric
problems thanks to its capability to drastically reduce the computational cost
required for the numerical
solutions~\cite{salmoiraghi2016advances,rozza2018advances}.
Despite the progresses in the global amount of computational power, many
problems still remain intractable using only the conventional discretization
methods --- e.g. finite element, finite volume --- especially in a many-queries
or real-time context.  The idea behind ROM is that a generic problem, even very
complex, has an intrinsic dimension much lower than the number of degrees of
freedom of the discretized system. To achieve this dimensionality reduction, a
database of several solution is firstly collected by solving the original
high-order model for different parameters (physical or geometrical).  Then, all
the solutions are combined to build the space onto which we can
accurately project the solution manifold and efficiently compute the solutions
for the new parameters. We typically call \emph{offline phase} the initial step in
which many high-fidelity solutions are computed and, depending from the studied
problem, it can be very demanding from the computational viewpoint. The second
step, the \emph{online phase}, is instead very fast since only the solution of the
low-dimensional problem has to be performed. 

The accuracy of the reduced order model depends by the problem itself, by the
number of parameters and the number of snapshots collected during the offline
phase. In this work, we use a non-intrusive proper orthogonal decomposition with
interpolation (PODI) method~\cite{bui2004aerodynamic} for the
data-driven dimensionality reduction, coupling it to the active
subspace (AS) property~\cite{constantine2015active}. 

PODI is an equation-free method based on proper orthogonal decomposition
capable to build a reduced order model without any knowledge about the
equations of the original problem, requiring just the $(\mupar,\mathbf{u})$
pairs, where $\mupar$ refers to the input parameters and $\mathbf{u}$
to the corresponding parametric solution. \TODO{We highlight that this method
does neither require information about the full-order formulation nor modifications
to the numerical solver}. Further insights and industrial applications of
data-driven \TODO{non-intrusive} ROMs can be found in~\cite{tezzele2018ecmi}.
Instead, AS is an emerging tool for the reduction of the parameter
space dimensionality. Basically, it aims to approximate a scalar function with
multi-dimensional input with a new function whose input parameters are linear
combinations of the original parameters. This parameter space reduction enhances
the accuracy of the reduced order model, thanks to the simplification of the
parametric formulation, overcoming the curse of dimensionality. 
For a more exhaustive introduction on ROMs and parameters space reduction, we
recommend~\cite{morhandbook2019}, while for the improvement of reduced order
models thanks to AS property, we mention~\cite{tezzele2019marine} for some
preliminary results in naval application using data-driven approach,
and~\cite{tezzele2018combined} as an example of enhanced POD-Galerkin method
applied to a biomedical problem.

In this work we are going to present an original application of the active
subspace property for the reconstruction of the modal coefficients of
the proper orthogonal decomposition. Moreover we demonstrate that this
results in better interpolation capabilities for the PODI when we have
a small amount of snapshots collected in the offline phase. This
coupling is especially useful when we have a very limited computational
budget.

The work is organised as follows: in Section~\ref{sec:pod} we present
the POD and the POD with interpolation as a data-driven approach for ROMs; in
Section~\ref{sec:active} we briefly introduce the active subspaces property; in
Section~\ref{sec:results} we show the application of the proposed
method on two simple parametric problems: a
structural analysis problem and a computational fluid dynamics
problem. Finally some perspective and future developments are presented.

\section{Reduced order modeling through proper orthogonal
  decomposition with interpolation}
\label{sec:pod}
\newcommand{\V}{\mathbb{V}}
\newcommand{\Ns}{{Ns}}
\newcommand{\uu}{\mathbf{u}}
\newcommand{\bfS}{\mathbf{S}}
\newcommand{\bfC}{\mathbf{C}}
\newcommand{\bfU}{\mathbf{U}}
\newcommand{\bfsigma}{\boldsymbol{\Sigma}}
\newcommand{\bfV}{\mathbf{V}}
\newcommand{\R}{\mathbb{R}}
\newcommand{\C}{\mathbb{A}}

Proper orthogonal decomposition (POD) is a widespread technique in the reduced
order modeling (ROM) community for the study of parametric problems, thanks to
its capacity to extract, from a set of high-dimensional snapshots, the basis
minimizing the error between the original snapshots and their orthogonal
projection. Let $\V^N$ be the high-dimensional discrete space which the
snapshots belong to. The basis $\{\phi_1, \phi_2, \dotsc, \phi_\Ns
\} \in \V^N$ spans the POD space $\V^\Ns \subset \V^N$ which, by construction,
is the optimal space of dimension $\Ns$ to represent the snapshots, where
$\phi_i$ for $i = 1, 2, \dotsc, \Ns$ are the so-called POD modes. 

We can project the equations of the full-order problem onto the POD space in
order to obtain a low-dimensional representation of the original operators.
Since the number of degrees of freedom $N$ of the initial problem is usually much
greater than the reduced dimension $\Ns$, the numerical solution of the reduced
order model results inexpensive from the computational viewpoint. This method
is called POD-Galerkin: for further details about this approach we suggest, among
many different works, for
example~\cite{hesthaven2016certified,ballarin2015supremizer,stabile2018finite,lorenzi2016pod,stabile2018reduced,couplet2005calibrated,karatzas2019reduced}. 

The data-driven approach \TODO{we present relies only on
  data and does not require the original equations describing the
  system. It is also non-intrusive in the sense that no modification
  of the simulation software is required (see for
example~\cite{zou2018nonintrusive} for a non-intrusive method that is
not data-driven)}. The original snapshots are projected onto the POD
space in order to reduce their dimensionality then the solution manifold is
approximated using an interpolation technique. Examples of application of this
so-called POD with interpolation (PODI)~\cite{bui2003proper,bui2004aerodynamic}
method can be found in literature: for naval engineering problem we
cite~\cite{demo2018isope,demo2018shape,demo2019marine}, for
automotive~\cite{salmoiraghi2018free,dolci2016proper} and for
aeronautics~\cite{ripepi2018reduced}. A coupling with isogeometric
analysis can be found in~\cite{garotta2018quiet}.

Let now focus on the computation of the POD modes. Let $\uu_i$, with
$i = 1,\dotsc, \Ns$, be the snapshots collected by solving the high-dimensional
problem, with different values of the input parameters $\mupar_i$, resulting in
$\Ns$ input-output pairs $(\mupar_i, \uu_i)$. The snapshots matrix $\bfS$ is
constructed arranging the snapshots as columns, such that $\bfS =
\begin{bmatrix}\uu_1 & \uu_2 & \dotsc & \uu_\Ns\\\end{bmatrix}$. We apply the
singular value decomposition \TODO{(see for
  example~\cite{brunton2019data})} to this matrix to obtain:
\begin{equation}
\bfS = \bfU \bfsigma \bfV^* \approx \bfU_k \bfsigma_k \bfV^*_k, 
\label{eq:svd}
\end{equation}
where $\bfU \in \C^{N\times\Ns}$ is the unitary matrix containing the
left-singular vectors, $\bfsigma \in \C^{\Ns\times\Ns}$ is the diagonal matrix
containing the singular values $\lambda_i$, and $\bfV \in \C^{\Ns\times\Ns}$,
with the symbol $^*$ denoting the conjugate transpose. The left-singular
vectors, namely the columns of $\bfU$, are the so-called POD modes. It is
important to note that the magnitude of the singular values describes the
energy of the corresponding modes. We can keep the first $k$ modes to span the
optimal space with dimension $k$ to represent the snapshots. Since the singular
values are returned in decreasing order, we can truncate the number of modes
simply selecting the first $k$ columns of $\bfU$.  The matrices $\bfU_k \in
\C^{N\times k}, \bfsigma_k \in \C^{k\times k}, \bfV_k \in \C^{\Ns\times k}$ in
Equation~\eqref{eq:svd} are the truncated matrices with rank $k$.  We can also
measure the error, in Euclidean and Frobenius norm, introduced by the
truncation~\cite{quarteroni2015reduced}:
\begin{equation}
\| \bfS - \bfU_k \bfsigma_k \bfV^*_k \|_2^2 = \lambda_{k+1}^2,
\end{equation}
\begin{equation}
\| \bfS - \bfU_k \bfsigma_k \bfV^*_k \|_F = \sqrt{\sum_{i=k+1}^\Ns{\lambda_i}^2}.
\end{equation}

After constructing the POD space, we can project the original snapshots onto this
space. In matrix form, we compute $\bfC \in
\R^{k\times \Ns}$ as $\bfC = \bfU_k^T \bfS$, where the columns of $\bfC$ are the low-dimensional
representation of the input, the so-called modal coefficients. Practically,
we can express the input snapshots as a linear combination of the modes using these coefficients.
Formally:
\begin{equation}
\uu_i = \sum_{j=1}^\Ns \alpha_{ji}\phi_j \approx \sum_{j=1}^k
\alpha_{ji}\phi_j,\quad\forall i \in [1, 2, \dotsc, \Ns],
\label{eq:lc}
\end{equation}
where $\alpha_{ji}$ are the elements of $\bfC$.
Finally, we obtain the $(\mupar_i, \alpha_i)$ pairs, for $i = 1, 2, \dotsc,
\Ns$, that sample the solution manifold in the parametric space. We are able
to interpolate the modal coefficients $\alpha$ and for any new parameter
approximate the new coefficients. At the end, we compute the high-dimensional solution by
projecting back the (approximated) modal coefficients to the original space by
using Equation~\eqref{eq:lc}.

Regarding the technical implementation of the PODI method, we adopt the Python package 
called EZyRB~\cite{demo2018ezyrb}.

\section{Active subspaces for modal coefficient reconstruction}
\label{sec:active}

Active subspaces
property~\cite{constantine2015active,constantine2014active} is a
powerful technique developed for parameter studies. The idea is to
perform a sensitivity analysis of the function of interest with
respect to the parameters, while reducing the parameter space
dimensionality.

We want to unveil specific directions in the parameter space along
which the scalar function of interest varies the most on
average. This is done by rescaling the inputs in a reference domain
centered in the origin and then by rotating the parameter space until
a lower dimensional structure is highlighted.
A drawback of this technique is that every function with a radial
symmetry does not have an active subspace since there are no preferred
directions to rotate the domain. Nevertheless a wide range of
functions of interest in engineering applications present an active
subspace of dimension one or two, resulting in a great reduction of
the parameter space. Among others we cite shape
optimization~\cite{lukaczyk2014active}, and uncertainty quantification
in the numerical simulation of the HyShot II
screamjet~\cite{constantine2015exploiting}. We would also like to mention
some naval engineering applications of the active subspaces: coupled
with boundary element method and free form
deformation~\cite{tezzele2018dimension}, for propeller blade
design~\cite{mola2019marine}, and a shared subspace application for
constraint optimisation~\cite{tezzele2018model}. 

Let $f$ be a parametric scalar function of interest $f(\mupar):
\mathbb{R}^p \to \mathbb{R}$, and $\rho: \mathbb{R}^p \to
\mathbb{R}^+$ a probability density function representing uncertainty
in the input parameters $\mupar \in \mathbb{R}^p$. Active subspaces
are a property of the pair $(f, \rho)$. To discover if a pair has an
active subspace of dimension $M$, we need to construct the uncentered
covariance matrix $\mathbf{\Sigma}$ of the gradients of $f$ with
respect to the input parameters, $\nabla f({\mupar}) \in
\mathbb{R}^p$, that reads:
\begin{equation}
\mathbf{\Sigma} = \mathbb{E}\, [\nabla_{\mupar} f \, \nabla_{\mupar} f
^T] =\int (\nabla_{\mupar} f) ( \nabla_{\mupar} f )^T
\rho \, d \mupar,
\end{equation}
where $\mathbb{E}$ is the expected value, and $\nabla_{\mupar} f
\equiv \nabla f({\mupar})$. Since $\mathbf{\Sigma}$ is symmetric it
has a real eigenvalue decomposition: 
\begin{equation}
\mathbf{\Sigma} = \mathbf{W} \Lambda \mathbf{W}^T,
\end{equation}
where $\mathbf{W}$ is an orthogonal matrix containing the eigenvectors
of $\mathbf{\Sigma}$ as columns, and $\Lambda$ is a diagonal matrix
composed by the non-negative eigenvalues arranged in descending order.

We are looking at spectral gap in order to identify the proper
dimension $M < p$ of the active subspace. In particular, we define the
active subspace of dimension $M$ as the span of the
first $M$ eigenvectors of $\mathbf{W}$, which correspond to the first
eigenvalues before a significan spectral gap. We proceed by
decomposing the two matrices as follows
\begin{equation}
\Lambda =   \begin{bmatrix} \Lambda_1 & \\
& \Lambda_2\end{bmatrix},
\qquad
\mathbf{W} = \left [ \mathbf{W}_1 \quad \mathbf{W}_2 \right ],
\qquad
\mathbf{W}_1 \in \mathbb{R}^{p\times M}.
\end{equation}
Then we can map the full parameters to the reduced ones through
$\mathbf{W}_1$. We define the active variable as $\mupar_M =
\mathbf{W}_1^T\mupar \in \mathbb{R}^M$, and the inactive variable as $\etapar =
\mathbf{W}_2^T\mupar \in \mathbb{R}^{p-M}$. The eigenvectors represent
the weights of the linear combination of the input parameters, thus
providing a sensitivity of each parameter. If a weight is almost zero,
that means $f$ does not vary along that direction on average. 

With the active variable we can build a surrogate model $g$ to approximate
the function of interest, that is
\begin{equation}
f (\mupar) \approx g(\mathbf{W}_1^T\mupar) = g (\mupar_M).
\end{equation}

The error upper bound for the approximation of $f$ through a response
surface depends on the square root of the sum of the eigenvalues
corresponding to the active variables times the eigenvectors
approximation error $\epsilon$, plus the square root of the sum of the remaining
eigenvalues~\cite{constantine2015active}:
\begin{equation}
\text{RMSE} \leq C_1 \left [ \epsilon \, \left(\sum_{i=1}^M \lambda_i \right)^{1/2} +
 \left ( \sum_{i=M+1}^p \lambda_i \right)^{1/2} \right ] + C_2,
\end{equation}
with $C_1$ and $C_2$ prescribed constants.

\section{Numerical results}
\label{sec:results}

In this section we are going to present two different test cases: the
first in the context of linear structural analysis, and the second in
computational fluid dynamics.

We sample the parameter space $\mathbb{P}$ with a uniform density
function, obtaining 200 samples for both the problems. The actual
space $\mathbb{P}$ will be defined in the corresponding sections.
To split the dataset for training and test purpose we use a $k$-fold
cross-validation (CV)~\cite{kohavi1995study}, with $k=5$. First we randomly
partition the samples into $k$ equal sized subsamples. Among these $k$ sets, a
single one is retained to validate the model, and the remaining $k-1$
are used as training data. We repeat the cross-validation process $k$ times, with
each of the $k$ subsamples used exactly once as test data. Then the
results are then averaged to produce a single estimation. One of the
advantages of this method over repeated random sub-sampling is that all
samples are used for both training and testing, and each observation
is used for validation exactly once.

We are interested in the relative reconstruction error of the two output fields
of interest: the stress tensor in the first example, and the pressure
field in the second one. The error is computed as the average over the
test samples of the norm of the difference between the exact and the
approximated solution over the norm of the exact solution.

We approximate the stress tensor field and the pressure field by using ridge
functions~\cite{pinkus1997approximating} $g_i$ to reconstruct the modal
coefficients through active subspaces.
The modal coefficients corresponding to the parameter sample
$\mupar^*$ are $\alpha(\mupar^*) \equiv \alpha^* = [\alpha_1^*, \dots,
\alpha_N^*]^T$.
The actual approximation of the modal coefficient is defined as follows
\begin{equation}
\label{eq:coef_approx}
\alpha_i^* \approx g_i(\mathbf{W}_{1, i}^T \, \mupar^*) \qquad i \in [1,
\dots, N],
\end{equation}
where $\mathbf{W}_{1, i}$ is the first eigenvector defining the active
subspace of dimension one corresponding to the $i$-th modal
coefficient, and $g_i$ a response function. We select only the first
eigenvector because every single modal coefficient presents an active
subspace of dimension one as shown in the following sections. 

Recalling Equation~\eqref{eq:lc}, we have the new approximated
representation of the snapshots $\uu_i$, using the corresponding $k$ response
functions $g_{ji}$:
\begin{equation}
\uu_i = \sum_{j=1}^\Ns \alpha_{ji}\phi_j \approx \sum_{j=1}^k
\alpha_{ji}\phi_j \approx \sum_{j=1}^k
g_{ji}(\mathbf{W}_{1, j}^T \, \mupar)\phi_j ,\quad\forall i \in [1, 2, \dotsc, \Ns].
\end{equation}

\subsection{Numerical study 1: stress tensor field reconstruction of a parametrised beam}
\begin{figure}[ht!]
\centering
\includegraphics[trim=30 200 30 300, clip, width=0.46\textwidth]{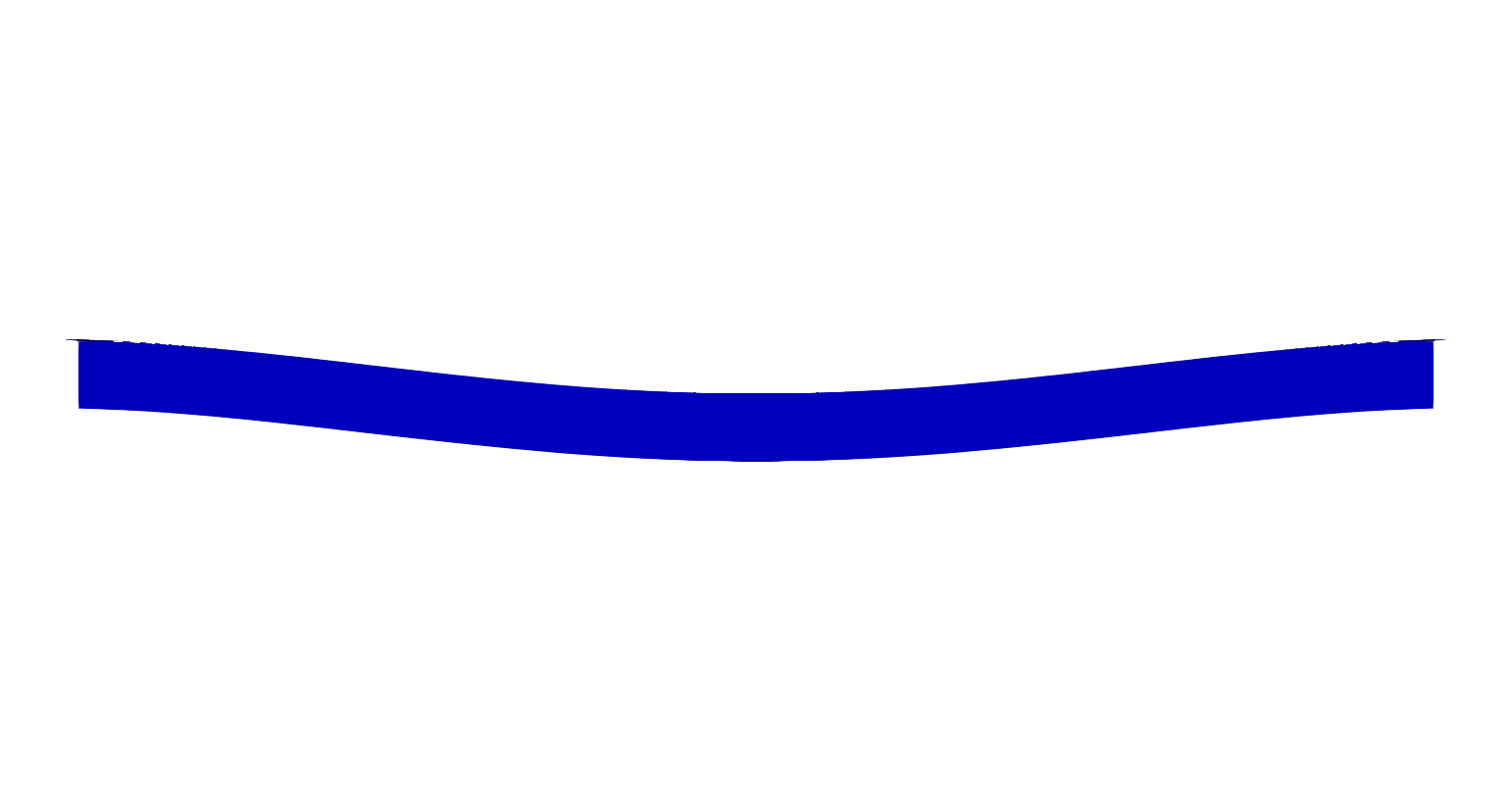}\hfill
\includegraphics[trim=40 150 80 150, clip, width=0.37\textwidth]{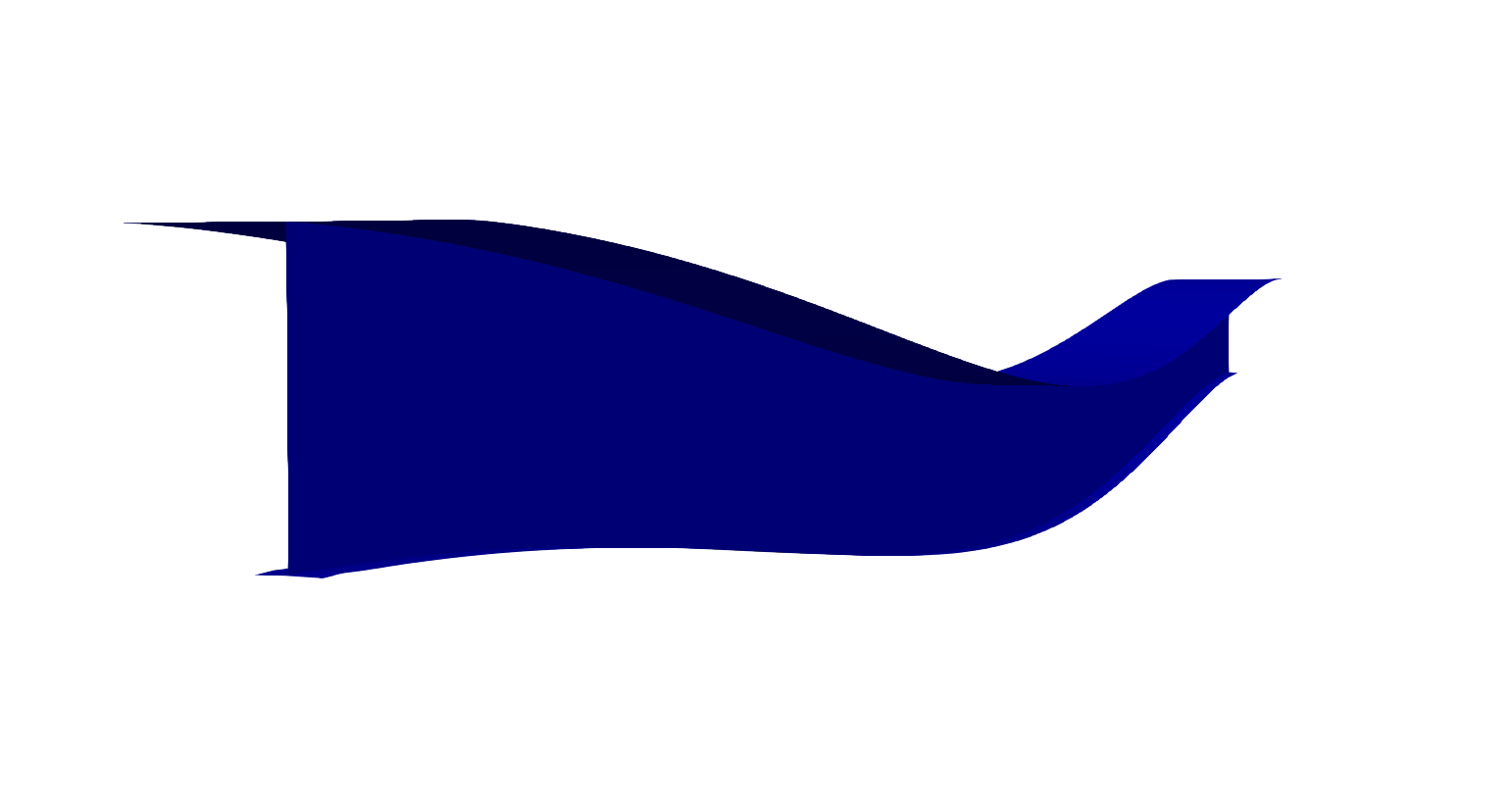}\hfill
  \begin{tikzpicture}
    \draw (0.1, 1.8) rectangle (2.4, 2.);
    \draw (1.15, 0.2) rectangle (1.35, 1.8);
    \draw (0.75, 0) rectangle (1.75, 0.2);
  \end{tikzpicture}
\caption{Two different views of a deformed beam after the application
  of a magnification factor equal to 100 along the z-direction. On the
right the longitudinal section of the beam.} 
\label{fig:def_beam}
\end{figure}
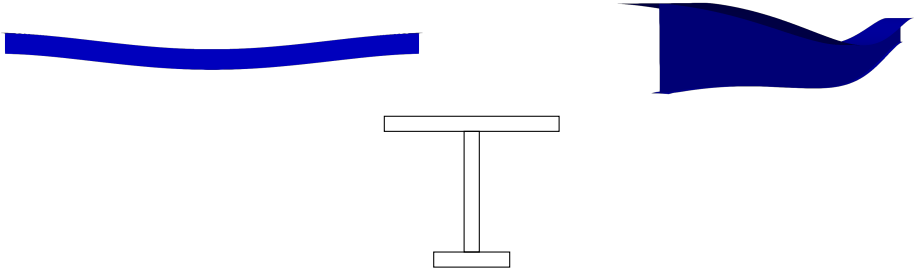

We consider a double-T beam with both ends fixed, and we
apply a uniform load condition. The length of the beam is 9000 mm, the
height of the beam is 450 mm, the span of the upper flange is equal to
500 mm, and the span of the lower flange is 100 mm. The uniform load
is 8.7 KN/m.
The input parameters $\mupar \in \mathbb{R}^p$ represent
the thickness of specific regions of the beam. We divided the beam in
three equidistant sections along the longitudinal direction, and we
preserve the symmetry by imposing the same thickness on the external
sections. The actual parameters are $\mupar \in \mathbb{P} := [5.0,
10.0]^6$, and the sampling is done using a uniform probability density
function. They are defined as in Table~\ref{tab:params}. The equations
governing the linear elastic isotropic problem are the equilibrium equation, the linearised
small-displacement strain-displacement relationship, and the Hooke's
law, respectively:
\begin{equation}
  \begin{cases}
 & - \nabla \cdot \sigma = f,\\
 & \epsilon = \frac{1}{2} [ \nabla u + \nabla u^T ], \\
 & \sigma = C(E, \nu) : \epsilon,
\end{cases}
\end{equation}
where $\sigma$ is the Cauchy stress tensor, $f$ is the body force,
$\epsilon$ is the infinitesimal strain tensor, $u$ is the displacement
vector, and $C$ is the fourth-order stiffness tensor depending on $E$,
the Young modulus, and on $\nu$, the Poisson's ratio.
In Figure~\ref{fig:def_beam} an example of a possible deformation of the
beam from two different views.

\begin{center}
\begin{table}[h]
\centering
\caption{Description of the input parameters $\mupar_i$, representing
  the thickness of particular sections of the beam. Lower and upper
  bounds are highlighted.}
\vspace{0.2cm}
\label{tab:params}
\begin{tabular}{cccc}
\toprule
  \qquad Parameter \qquad \qquad & Section of the beam  & \quad Lower bound (mm) \quad \quad  & \quad Upper bound (mm) \quad \quad \\ \midrule
$\mu_1$ & External sections web & 5.0 & 10.0 \\
$\mu_2$ & Internal section web &  5.0 & 10.0 \\
$\mu_3$ & External sections lower flange & 5.0 & 10.0 \\
$\mu_4$ & Internal section lower flange &  5.0 & 10.0 \\ %\midrule
$\mu_5$ & External sections upper flange &  5.0 & 10.0 \\
$\mu_6$ & Internal section upper flange &  5.0 & 10.0 \\ \bottomrule
\end{tabular}
\end{table}
\end{center}

\begin{figure}[h!]
\centering
\includegraphics[width=0.485\textwidth]{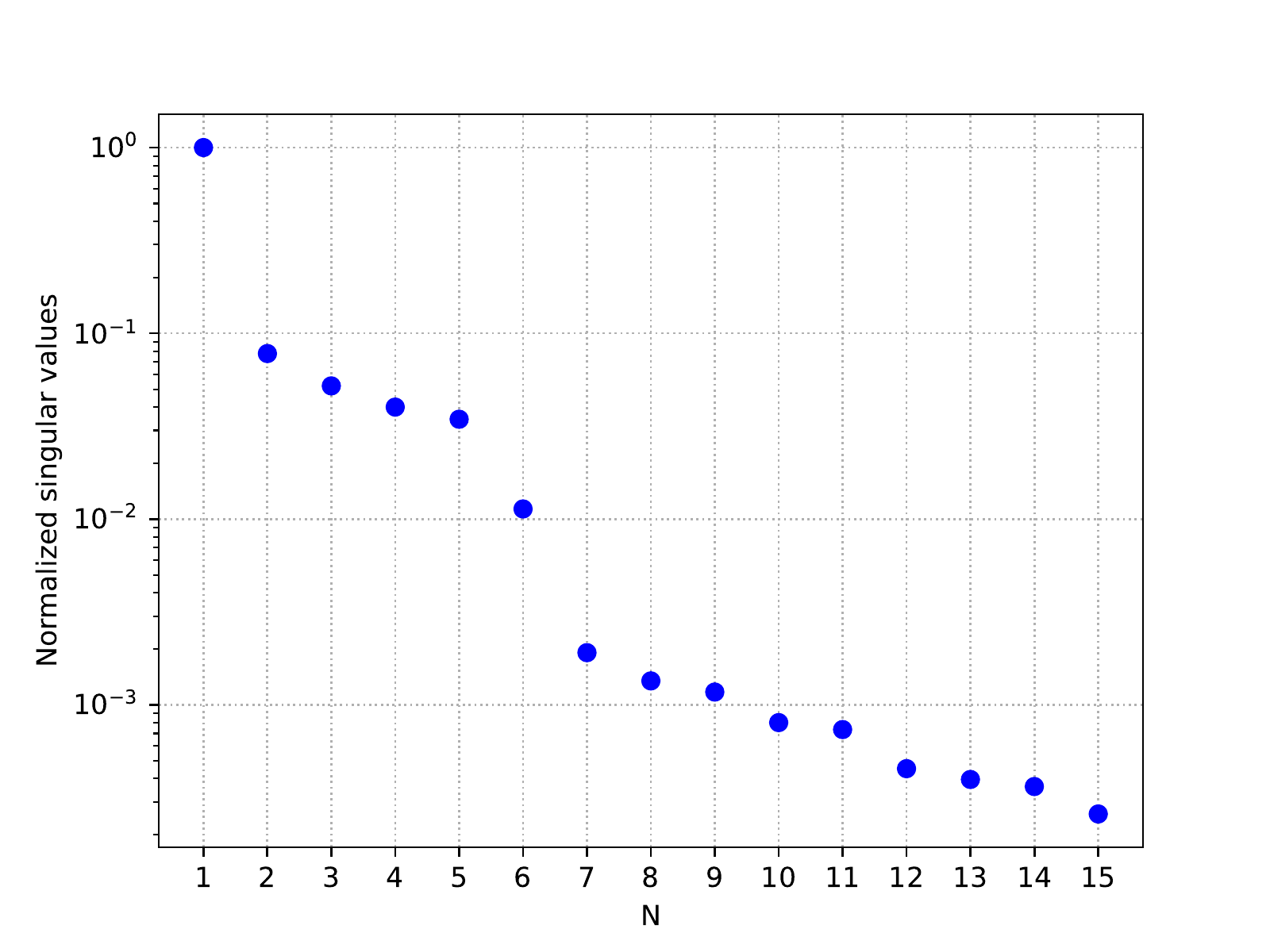}\hfill
\includegraphics[width=0.485\textwidth]{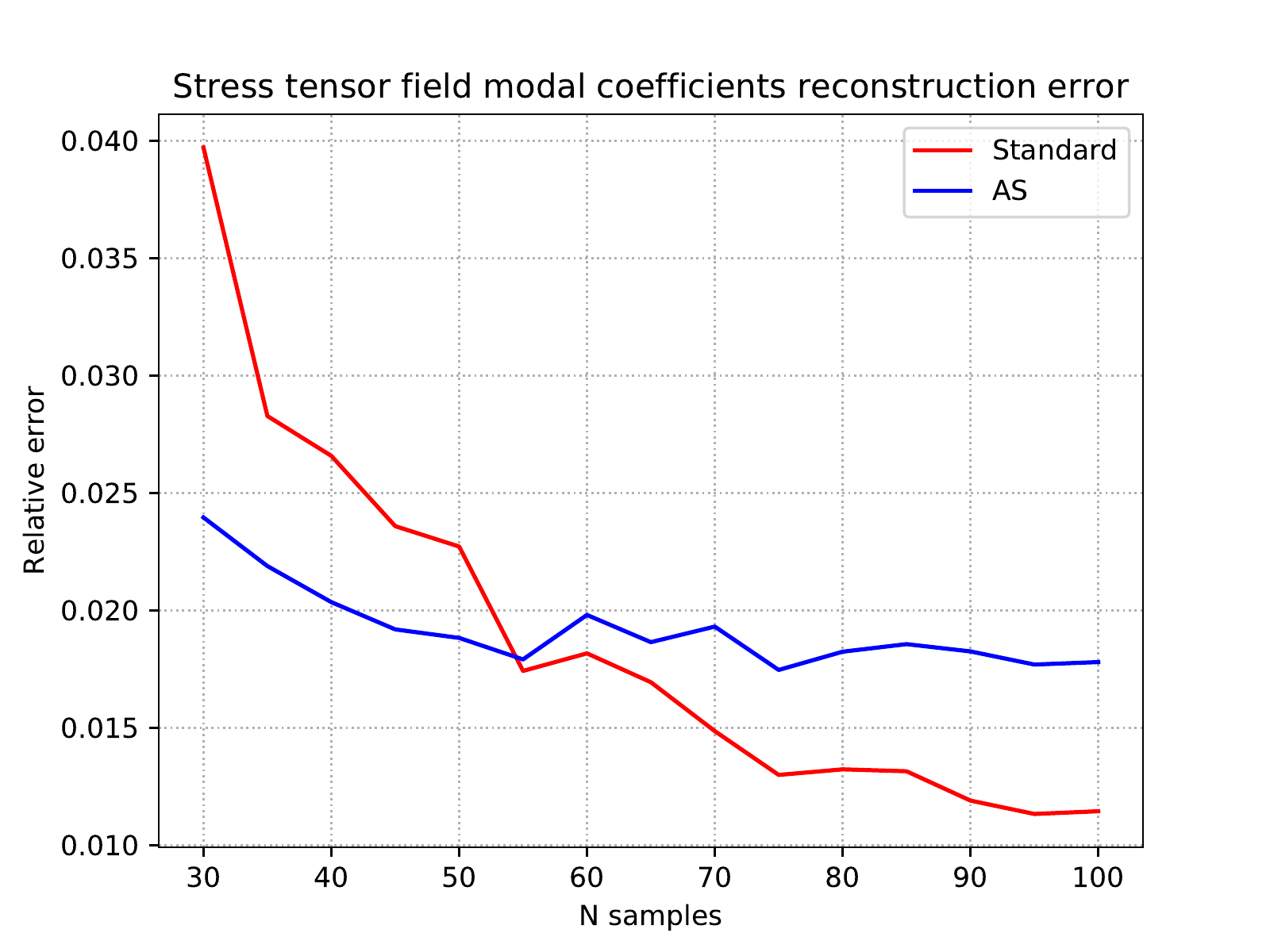}
\caption{On the left the normalized singular values decay of the
  stress tensor field snapshots matrix. Only the first 15 are depicted, and we select only
  the first 6. On the right the relative error for the stress tensor field
  modal coefficients reconstruction using 5-fold cross validation for
  both the full parameter space and the active subspace (AS), as a
  function of the number of the samples.} 
\label{fig:stresses}
\end{figure}

We compute the POD modes and coefficients as explained in
Section~\ref{sec:pod}. In particular we truncate the modal expansion
to the first 6 modes. This choice is made in accordance to the singular
values decay as depicted in Figure~\ref{fig:stresses} on the left. After the
computation of the modal coefficients we approximate the active
subspace of dimension 1 for each coefficient. Following
Equation~\eqref{eq:coef_approx} we are able to approximate the
coefficients using a univariate function. As response surface $g$ we use
gaussian process regression with a radial basis function kernel. The
comparison with the interpolation of the single coefficients using the
full parameter space is presented in
Figure~\ref{fig:stresses} on the right. For the \emph{standard} procedure
we mean the use of radial basis functions interpolation in order to
approximate the functions $\mupar \in \mathbb{P} \to
\alpha_i(\mupar)$, with $i \in [1, \dots, 6]$. The relative error is
computed as describe above using a 5-fold cross validation method, and
is expressed as a function of the cardinality of the samples set (both
for training and testing). We can see that active subspaces allow
smaller reconstruction error for a given number of samples up to a
certain threshold, 44 training samples in this case. This means that,
for a given small computational budget devoted to the offline phase, we can
gain around 1\% in terms of relative error. The reconstruction error
in the AS case does not improve increasing the samples because we
have already reached intrinsic error of the active subspaces
approximation.
% \begin{figure}
% \centering
% \includegraphics[width=0.5\textwidth]{figures/plot_as_gpy.pdf}
% \caption{The relative error for the stress tensor field
%   modal coefficients reconstruction using 5-fold cross validation for
%   both the full parameter space and the active subspace (AS), as a
%   function of the number of the samples.} 
% \label{fig:interp_err_stress}
% \end{figure}

\subsection{Numerical study 2: pressure reconstruction in fluid
  dynamics with a parametric domain}

We are now considering the flow around a square in a parametric domain, a
standard benchmark in computational fluid dynamics. Since the used approach is
completely data-driven, it is out of the purpose of this work the analysis of
the mathematical model. We just provide an overview of the
parametric formulation to give the opportunity to the reader to understand
and replicate the results.

We have a rectangular domain in which an incompressible flow past a square
cylinder. A sketch of the domain is depicted in Figure~\ref{fig:domain}. The
$\Gamma_\text{in}$, $\Gamma_\text{out}$ and $\Gamma_\text{wall}$ refer
respectively to the inlet boundary, the outlet boundary and the physical wall
of the domain. To stay in a steady regime, we impose a low Reynolds number ($Re
= 13$). The incompressible Navier Stokes equations
\begin{equation}
\begin{cases}
\frac{\partial u}{\partial t} + (u \cdot \nabla)u - \nabla \cdot
\nu\nabla u = - \nabla p & \text{on}\, \Omega ,\\
\nabla \cdot u = 0& \text{on}\, \Omega ,\\
\end{cases}
\end{equation}
are solved in the described domain using OpenFOAM~\cite{of}, a open source library
implementing finite volume (FV) method.

\pgfplotsset{
    standard/.style={
        axis x line=middle,
        axis y line=left,
        ticklabel style={fill=white},
        set layers=tick labels on top% use layers and choose the new layer set
    },
    layers/tick labels on top/.define layer set=% define the new layer set based on the standard one
        {axis background,axis grid,axis ticks,axis lines,main,%
          axis tick labels,% <- tick labels before main
          axis descriptions,axis foreground}
        {/pgfplots/layers/standard}
}

\begin{figure}[h!]
  \centering
  \begin{tikzpicture}%[scale=75]

    \begin{axis}[
        ymax =  0.05,
        ymin = -0.05,
        xmax =  0.15,
        xmin = -0.05,
        unit vector ratio*=1 1 1,
        width=13cm,
        grid=major,
        grid style={dashed, gray!30},
        scaled x ticks = false,
        legend pos=south east,
        xticklabel style={/pgf/number format/fixed}
    ]
    \addplot[black, -, line width=.2mm, forget plot] coordinates {
        (-0.005,   -0.005)
        ( 0.005,   -0.005)
        ( 0.005,    0.005)
        (-0.005,    0.005)
        (-0.005,   -0.005)
    };
    \begin{scope}
    \foreach \x in {-0.0350000000000000000,-0.0316666666666666667,-0.0283333333333333334,-0.0250000000000000001,-0.0216666666666666668,-0.0183333333333333335,-0.0150000000000000002,-0.0116666666666666669,-0.0083333333333333336,-0.0050000000000000003,-0.0016666666666666670,0.0016666666666666663,0.0049999999999999996,0.0083333333333333329,0.0116666666666666662,0.0149999999999999995,0.0183333333333333328,0.0216666666666666661,0.0249999999999999994,0.0283333333333333327,0.0316666666666666660,0.0349999999999999993,0.0383333333333333326,0.0416666666666666659,0.0449999999999999992,0.0483333333333333325,0.0516666666666666658,0.0549999999999999991,0.0583333333333333324,0.0616666666666666657,0.0649999999999999990,0.0683333333333333323,0.0716666666666666656,0.0749999999999999989,0.0783333333333333322,0.0816666666666666655,0.0849999999999999988,0.0883333333333333321,0.0916666666666666654,0.0949999999999999987} {
    \foreach \y in {-0.0350000000000000000,-0.0316666666666666667,-0.0283333333333333334,-0.0250000000000000001,-0.0216666666666666668,-0.0183333333333333335,-0.0150000000000000002,-0.0116666666666666669,-0.0083333333333333336,-0.0050000000000000003,-0.0016666666666666670,0.0016666666666666663,0.0049999999999999996,0.0083333333333333329,0.0116666666666666662,0.0149999999999999995,0.0183333333333333328,0.0216666666666666661,0.0249999999999999994,0.0283333333333333327,0.0316666666666666660,0.034999999999999999} 
    {
    %\node[circle] (pts) at (axis cs:\x,\y) {};
    %\fill[black!30] (axis cs:{\x},\y) circle[radius=0.5pt];
    \addplot [only marks,mark=o, mark size=0.5, fill=white, gray, forget plot] coordinates { (\x,\y) }; 
    }}
    \end{scope}

    \addplot [only marks,mark=o, mark size=0.5, fill=white, gray] coordinates { (-0.035, -0.035) }; 
    \addlegendentry{FFD candidate control points}

    \addplot [only marks,mark=o, mark size=1.5, fill=white, red, forget plot] coordinates { (-0.005, -0.005) }; 
    \addplot [only marks,mark=o, mark size=1.5, fill=white, red, forget plot] coordinates { (0.005, -0.005) }; 
    \addplot [only marks,mark=o, mark size=1.5, fill=white, red, forget plot] coordinates { (0.005, 0.005) }; 
    \addplot [only marks,mark=o, mark size=1.5, fill=white, red] coordinates { (-0.005, 0.005) }; 
    \addlegendentry{moved FFD control points}

    \node (v1) at (axis cs:-0.035, -0.035) {};
    \node (v2) at (axis cs: 0.095, -0.035) {};
    \node (v3) at (axis cs: 0.095,  0.035) {};
    \node (v4) at (axis cs:-0.035,  0.035) {};

    \node (v5) at (axis cs:-0.005, -0.005) {};
    \node (v6) at (axis cs: 0.005, -0.005) {};
    \node (v7) at (axis cs: 0.005,  0.005) {};
    \node (v8) at (axis cs:-0.005,  0.005) {};

    \node (v5x) at (axis cs:-0.012, -0.005) {};
    \node (v6x) at (axis cs: 0.012, -0.005) {};
    \node (v7x) at (axis cs: 0.012,  0.005) {};
    \node (v8x) at (axis cs:-0.012,  0.005) {};

    \node (v5y) at (axis cs:-0.005, -0.012) {};
    \node (v6y) at (axis cs: 0.005, -0.012) {};
    \node (v7y) at (axis cs: 0.005,  0.012) {};
    \node (v8y) at (axis cs:-0.005,  0.012) {};

    \node (omega) at (axis cs:0.09, 0.03) {$\Omega$};

    \draw[-, line width=0.2mm] (v1.center) -- (v2.center) node[midway, above] {$\Gamma_\text{wall}$};
    \draw[-, line width=0.2mm] (v1.center) -- (v4.center) node[midway, left] {$\Gamma_\text{in}$};
    \draw[-, line width=0.2mm] (v4.center) -- (v3.center) node[midway, above] {$\Gamma_\text{wall}$};
    \draw[-, line width=0.2mm] (v2.center) -- (v3.center) node[midway, left] {$\Gamma_\text{out}$};

    \draw[->, red, line width=0.3mm] (v5.center) -- (v5x.center) node[fill=white, left] {$\mupar_1$};
    \draw[->, red, line width=0.3mm] (v6.center) -- (v6x.center) node[right] {$\mupar_2$};
    \draw[->, red, line width=0.3mm] (v7.center) -- (v7x.center) node[right] {$\mupar_3$};
    \draw[->, red, line width=0.3mm] (v8.center) -- (v8x.center) node[left] {$\mupar_4$};

    \draw[->, red, line width=0.3mm] (v5.center) -- (v5y.center) node[below] {$\mupar_5$};
    \draw[->, red, line width=0.3mm] (v6.center) -- (v6y.center) node[below] {$\mupar_6$};
    \draw[->, red, line width=0.3mm] (v7.center) -- (v7y.center) node[above] {$\mupar_7$};
    \draw[->, red, line width=0.3mm] (v8.center) -- (v8y.center) node[above] {$\mupar_8$};

    \end{axis}
    \end{tikzpicture}
\caption{A sketch of the (undeformed) domain for the computational fluid
dynamics simulation. The black dots are the FFD candidate control points, while
in red we highlight the ones we move to morph the original domain.}
\label{fig:domain}
\end{figure}
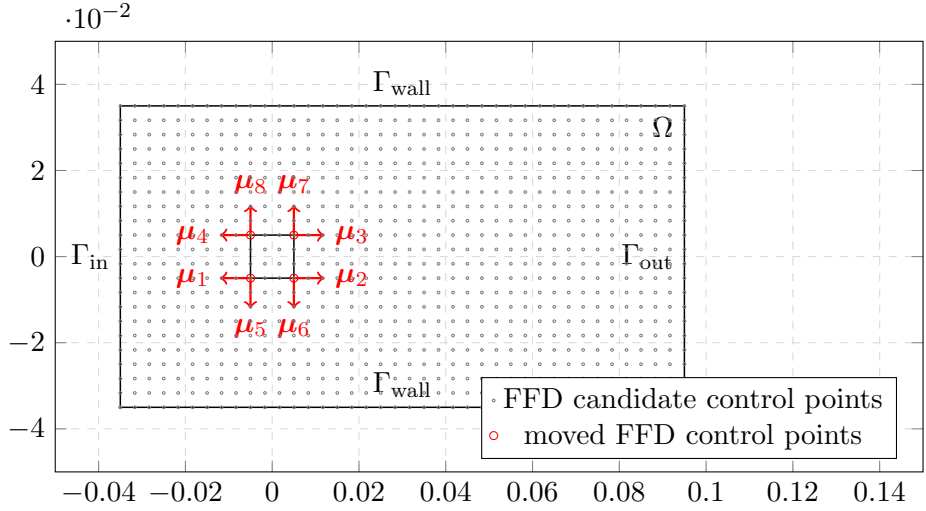

To parametrise and deform the original domain, we apply the free form deformation (FFD)
technique~\cite{sederbergparry1986} to the undeformed computational grid.
It is a well-spread parameterisation method used to morph smoothly complex geometries.
It consists on the composition of three maps:
\begin{enumerate}
\item a function $\psi$ that trasform the physical domain $\Omega$ to
  the reference domain $\widehat{\Omega}$;
\item the actual morphing done by the map $T$. It induces the
  trasformation of $\widehat{\Omega}$ given the displacements of the
  FFD control points composing the lattice around the object to
  deform. It is based on tensor products of Bernstein
  polynomials. This produces  $\widehat{\Omega} (\mupar)$; 
\item the back mapping from the reference domain to the physical one:
  $\psi^{-1}: \widehat{\Omega} (\mupar) \to \Omega (\mupar)$.
\end{enumerate}
The actual mapping takes the following form: $\mathcal{M} (\mupar): \psi^{-1}
\circ T \circ \psi \,(\mupar)$.

\begin{figure}[b!]
\centering
\includegraphics[width=.485\textwidth]{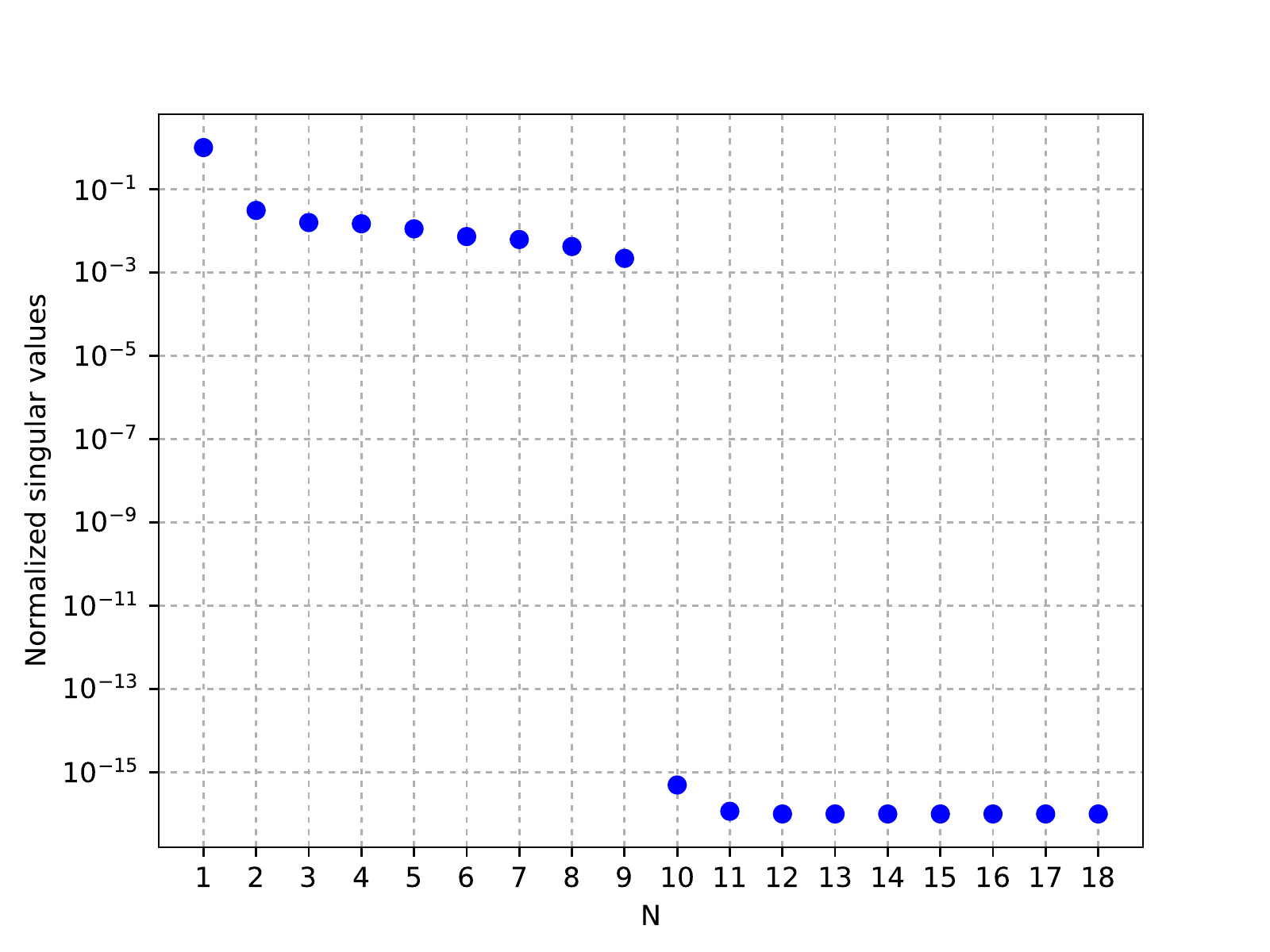}\hfill
\includegraphics[width=.485\textwidth]{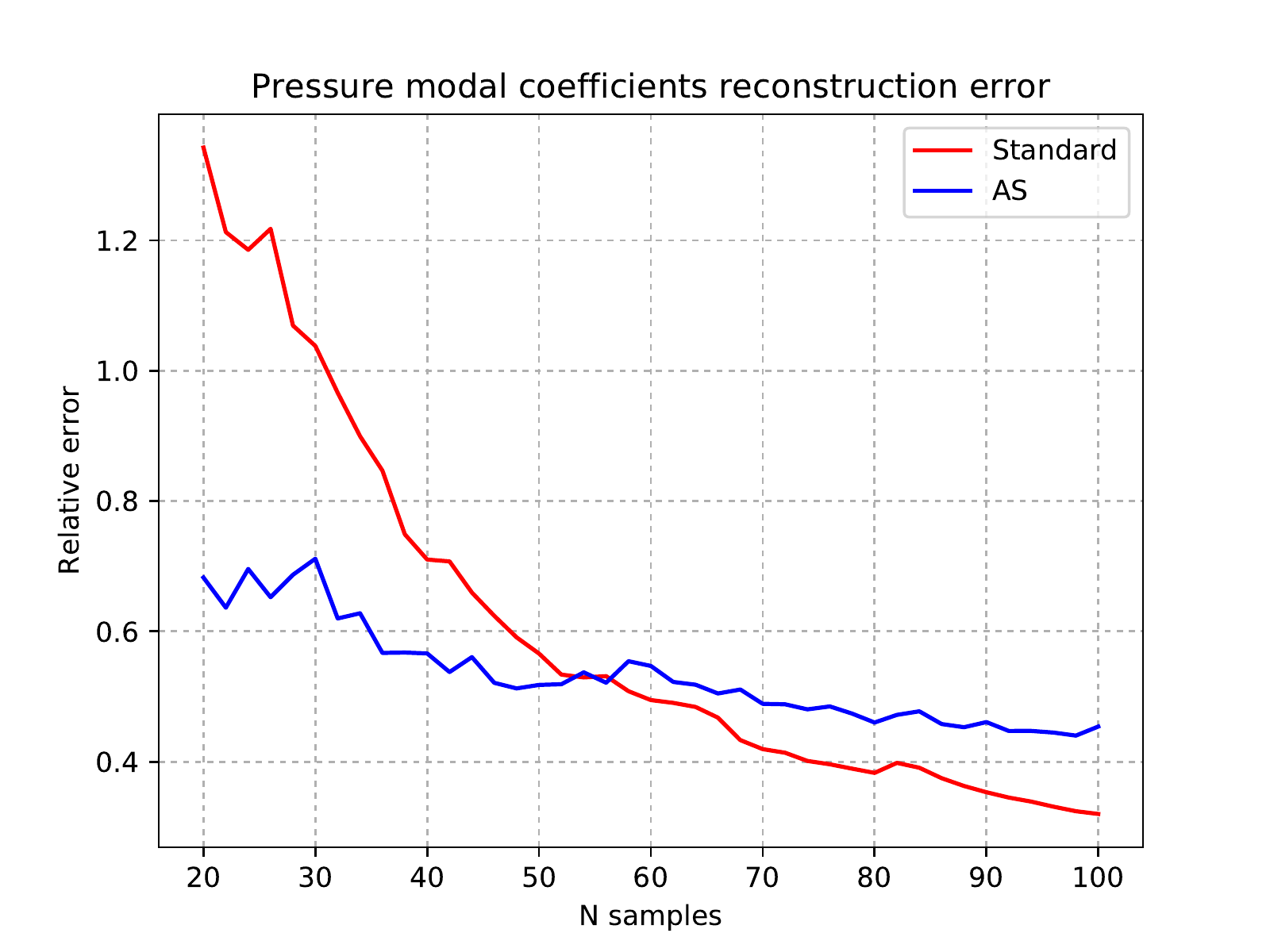}
\caption{On the left the normalized singular values decay of the
  pressure snapshots. On the right the relative error for the pressure modal coefficients
  reconstruction using 5-fold cross validation for both the full
  parameter space and the active subspace (AS), as a function of the
  number of the samples.} 
\label{fig:pressure}
\end{figure}

In this work, we embedded all the domain into the FFD lattice of control
points, moving only the 4 control points collocated on the internal square
vertices. The number of parameters $p$ is set to be equal to 8. This is the
result of moving 4 FFD control points along the $x$- and $y$-directions (in
Figure~\ref{fig:domain} the direction of the displacements is reported). The
computational grid is deformed using the Python package PyGeM~\cite{pygem},
implementing some of the most popular morphing technique as FFD, radial basis
functions (RBF) interpolation, and inverse distance weighting (IDW). For more
details about these methodologies, we
suggest~\cite{buhmann2003radial,forti2014efficient,rozza2013free,shepard1968,ballarin2019pod}.

\begin{figure}[h!]
\centering
\includegraphics[trim=120 0 120 0, clip, width=1.\textwidth]{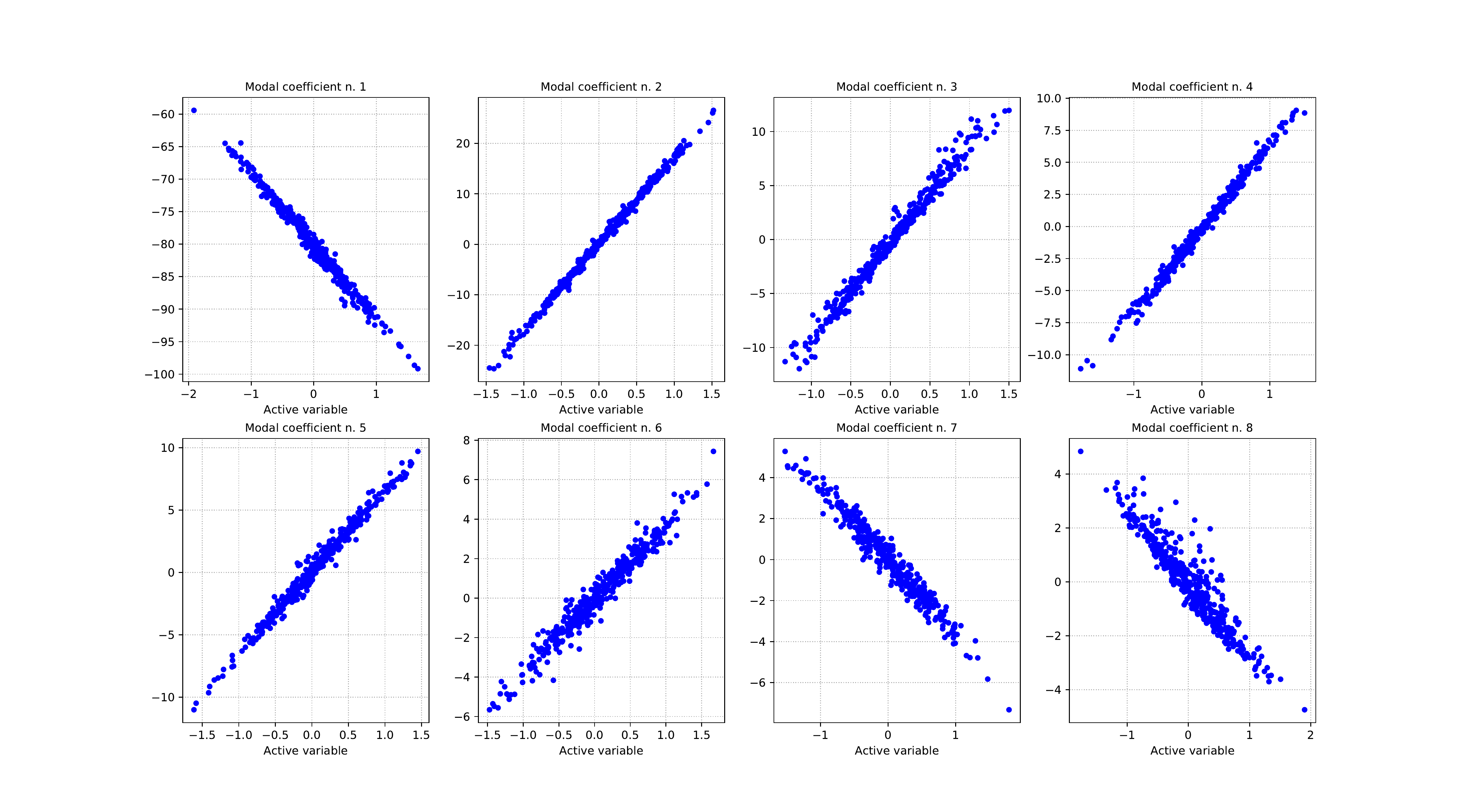}
\caption{Sufficient summary plots for the pressure modal
  coefficients. The first 8 coefficients $\alpha_i (\mupar)$ are represented along the
  corresponding active variable $\mathbf{W}_{1, i}^T \, \mupar$.} 
\label{fig:ssp_pressure}
\end{figure}
Once we collected all the numerical solutions for the 200 sampling
configurations, we extract the pressure snapshots to compute the POD modes.
Analyzing the singular values decay (Figure~\ref{fig:pressure} on the left), we select the
first \TODO{9} modes. After the projection of the samples onto the POD space, the
active subspace of dimension 1 is then computed for each modal coefficient.
Figure~\ref{fig:ssp_pressure} reports the AS accuracy by showing the
approximation of the modal coefficients along the new active
variable, that is $\alpha_i (\mupar)$ against $\mathbf{W}_{1, i}^T \,
\mupar$. As in the previous example, a radial basis function kernel
has been chosen to perform the gaussian process regression to reconstruct the response
surface in the active subspaces.  We compare the accuracy of the data-driven
model with the AS enhancement and using a standard approach: the latter refers
to an interpolation in the full dimensional parametric space ($\R^8$) of the
modal coefficients using a radial basis function. We compute the relative error
with a CV technique, varying the number of snapshots used for the AS
approximation and for POD modes computation (Figure~\ref{fig:pressure}
on the right). Using AS, the relative error is lower with respect to the
standard approach with a limited set of high-fidelity
snapshots and with only 20 samples the relative error is halved thanks to AS.
Increasing the initial database dimension, the accuracy gain becomes less
important and the standard interpolation provides better results, in this
example using more than 55 initial samples, since the error introduced by AS
approximation is greater than the interpolation error. 
% \begin{figure}[h!]
% \centering
% \includegraphics[width=.5\textwidth]{figures/error_p.pdf}
% \caption{The relative error for the pressure modal coefficients
%   reconstruction using 5-fold cross validation for both the full
%   parameter space and the active subspace (AS), as a function of the
%   number of the samples.}
% \label{fig:error_p}
% \end{figure}

\section{Conclusions and perspectives}
\label{sec:end}

In this work we presented a novel coupling between the active subspaces
property and the non-intrusive reduced order method called proper
orthogonal decomposition with interpolation. We proved the efficiency
of the technique on two diverse benchmark problems involving
geometrical parameterisation. The AS property is able to enhance the
performance of PODI when the offline database of solutions is poor in
terms of cardinality. This results in better reconstruction error of
the output fields of interest for new untried input parameters. So the
proposed method is useful especially when the offline phase of the
reduced order model is computationally intensive. Nevertheless after a
certain threshold, that means having enough samples of the full order
solutions, the classical PODI over the full parameter space is
more viable.

Future developments are required in order to have error bounds for the
proposed method. Moreover, non linear extensions of the active
subspaces property could be beneficial to improve the approximation
accuracy.

\section*{Acknowledgements}
This work was partially supported by an industrial Ph.D. grant sponsored
by Fincantieri S.p.A., \TODO{by the INdAM-GNCS 2019 project ``Advanced
intrusive and non-intrusive model order reduction techniques and
applications'',} and by the European Union Funding for Research and
Innovation --- Horizon 2020 Program --- in the framework
of European Research Council Executive Agency: H2020 ERC CoG 2015
AROMA-CFD project 681447 ``Advanced Reduced Order Methods with
Applications in Computational Fluid Dynamics'' P.I. Gianluigi
Rozza.

\end{document}